\documentclass[10pt]{amsart}

\makeatletter
 
\def\diagram{\leftwidth=\z@ \rightwidth=\z@ \topheight=\z@
\botheight=\z@ \setbox\@picbox\hbox\bgroup}
 
\def\enddiagram{\egroup\wd\@picbox\rightwidth\unitlength
\ht\@picbox\topheight\unitlength \dp\@picbox\botheight\unitlength
\hskip\leftwidth\unitlength\box\@picbox}
 
\def\bfig{\begin{diagram}}
\def\efig{\end{diagram}}
\newcount\wideness \newcount\leftwidth \newcount\rightwidth
\newcount\highness \newcount\topheight \newcount\botheight
 
\def\ratchet#1#2{\ifnum#1<#2 \global #1=#2 \fi}
 
\def\putbox(#1,#2)#3{%
\horsize{\wideness}{#3} \divide\wideness by 2
{\advance\wideness by #1 \ratchet{\rightwidth}{\wideness}}
{\advance\wideness by -#1 \ratchet{\leftwidth}{\wideness}}
\vertsize{\highness}{#3} \divide\highness by 2
{\advance\highness by #2 \ratchet{\topheight}{\highness}}
{\advance\highness by -#2 \ratchet{\botheight}{\highness}}
\put(#1,#2){\makebox(0,0){$#3$}}}
 
\def\putlbox(#1,#2)#3{%
\horsize{\wideness}{#3}
{\advance\wideness by #1 \ratchet{\rightwidth}{\wideness}}
{\ratchet{\leftwidth}{-#1}}
\vertsize{\highness}{#3} \divide\highness by 2
{\advance\highness by #2 \ratchet{\topheight}{\highness}}
{\advance\highness by -#2 \ratchet{\botheight}{\highness}}
\put(#1,#2){\makebox(0,0)[l]{$#3$}}}
 
\def\putrbox(#1,#2)#3{%
\horsize{\wideness}{#3}
{\ratchet{\rightwidth}{#1}}
{\advance\wideness by -#1 \ratchet{\leftwidth}{\wideness}}
\vertsize{\highness}{#3} \divide\highness by 2
{\advance\highness by #2 \ratchet{\topheight}{\highness}}
{\advance\highness by -#2 \ratchet{\botheight}{\highness}}
\put(#1,#2){\makebox(0,0)[r]{$#3$}}}

\def\adjust[#1]{} 
 
\newcount \coefa
\newcount \coefb
\newcount \coefc
\newcount\tempcounta
\newcount\tempcountb
\newcount\tempcountc
\newcount\tempcountd
\newcount\xext
\newcount\yext
\newcount\xoff
\newcount\yoff
\newcount\gap%
\newcount\arrowtypea
\newcount\arrowtypeb
\newcount\arrowtypec
\newcount\arrowtyped
\newcount\arrowtypee
\newcount\height
\newcount\width
\newcount\xpos
\newcount\ypos
\newcount\run
\newcount\rise
\newcount\arrowlength
\newcount\halflength
\newcount\arrowtype
\newdimen\tempdimen
\newdimen\xlen
\newdimen\ylen
\newsavebox{\tempboxa}%
\newsavebox{\tempboxb}%
\newsavebox{\tempboxc}%
 
\newdimen\w@dth
 
\def\setw@dth#1#2{\setbox\z@\hbox{$#1$}\w@dth=\wd\z@
\setbox\@ne\hbox{$#2$}\ifnum\w@dth<\wd\@ne \w@dth=\wd\@ne \fi
\advance\w@dth by 1.2em}
 
\def\t@^#1_#2{\def\n@one{#1}\def\n@two{#2}\mathrel{\setw@dth{#1}{#2}
\mathop{\hbox to \w@dth{\rightarrowfill}}\limits
\ifx\n@one\empty\else ^{\box\z@}\fi
\ifx\n@two\empty\else _{\box\@ne}\fi}}
\def\t@@^#1{\@ifnextchar_ {\t@^{#1}}{\t@^{#1}_{}}}
\def\to{\@ifnextchar^ {\t@@}{\t@@^{}}}
 
\def\t@left^#1_#2{\def\n@one{#1}\def\n@two{#2}\mathrel{\setw@dth{#1}{#2}
\mathop{\hbox to \w@dth{\leftarrowfill}}\limits
\ifx\n@one\empty\else ^{\box\z@}\fi
\ifx\n@two\empty\else _{\box\@ne}\fi}}
\def\t@@left^#1{\@ifnextchar_ {\t@left^{#1}}{\t@left^{#1}_{}}}
\def\toleft{\@ifnextchar^ {\t@@left}{\t@@left^{}}}
 
\def\two@^#1_#2{\def\n@one{#1}\def\n@two{#2}\mathrel{\setw@dth{#1}{#2}
\mathop{\vcenter{\hbox to \w@dth{\rightarrowfill}\kern-1.7ex
                 \hbox to \w@dth{\rightarrowfill}}%
       }\limits
\ifx\n@one\empty\else ^{\box\z@}\fi
\ifx\n@two\empty\else _{\box\@ne}\fi}}
\def\tw@@^#1{\@ifnextchar_ {\two@^{#1}}{\two@^{#1}_{}}}
\def\two{\@ifnextchar^ {\tw@@}{\tw@@^{}}}
 
\def\tofr@^#1_#2{\def\n@one{#1}\def\n@two{#2}\mathrel{\setw@dth{#1}{#2}
\mathop{\vcenter{\hbox to \w@dth{\rightarrowfill}\kern-1.7ex
                 \hbox to \w@dth{\leftarrowfill}}%
       }\limits
\ifx\n@one\empty\else ^{\box\z@}\fi
\ifx\n@two\empty\else _{\box\@ne}\fi}}
\def\t@fr@^#1{\@ifnextchar_ {\tofr@^{#1}}{\tofr@^{#1}_{}}}
\def\tofro{\@ifnextchar^ {\t@fr@}{\t@fr@^{}}}

\def\mon{\mathop{\m@th\hbox to
      14.6\P@{\lasyb\char'51\hskip-2.1\P@$\arrext$\hss
$\mathord\rightarrow$}}\limits} 
\def\leftmono{\mathrel{\m@th\hbox to
14.6\P@{$\mathord\leftarrow$\hss$\arrext$\hskip-2.1\P@\lasyb\char'50%
}}\limits} 
\mathchardef\arrext="0200       

\def\settypes(#1,#2,#3){\arrowtypea#1 \arrowtypeb#2 \arrowtypec#3}
\def\settoheight#1#2{\setbox\@tempboxa\hbox{#2}#1\ht\@tempboxa\relax}%
\def\settodepth#1#2{\setbox\@tempboxa\hbox{#2}#1\dp\@tempboxa\relax}%
\def\settokens[#1`#2`#3`#4]{%
     \def\tokena{#1}\def\tokenb{#2}\def\tokenc{#3}\def\tokend{#4}}
\def\setsqparms[#1`#2`#3`#4;#5`#6]{%
\arrowtypea #1
\arrowtypeb #2
\arrowtypec #3
\arrowtyped #4
\width #5
\height #6
}
\def\setpos(#1,#2){\xpos=#1 \ypos#2}

\def\settriparms[#1`#2`#3;#4]{\settripairparms[#1`#2`#3`1`1;#4]}%
 
\def\settripairparms[#1`#2`#3`#4`#5;#6]{%
\arrowtypea #1
\arrowtypeb #2
\arrowtypec #3
\arrowtyped #4
\arrowtypee #5
\width #6
\height #6
}
 
\def\resetparms{\settripairparms[1`1`1`1`1;500]\width 500}
 
\resetparms
 
\def\mvector(#1,#2)#3{
\put(0,0){\vector(#1,#2){#3}}%
\put(0,0){\vector(#1,#2){26}}%
}
\def\evector(#1,#2)#3{{
\arrowlength #3
\put(0,0){\vector(#1,#2){\arrowlength}}%
\advance \arrowlength by-30
\put(0,0){\vector(#1,#2){\arrowlength}}%
}}
 
\def\horsize#1#2{%
\settowidth{\tempdimen}{$#2$}%
#1=\tempdimen
\divide #1 by\unitlength
}
 
\def\vertsize#1#2{%
\settoheight{\tempdimen}{$#2$}%
#1=\tempdimen
\settodepth{\tempdimen}{$#2$}%
\advance #1 by\tempdimen
\divide #1 by\unitlength
}
 
\def\putvector(#1,#2)(#3,#4)#5#6{{%
\ifnum3<\arrowtype
\putdashvector(#1,#2)(#3,#4)#5\arrowtype
\else
\ifnum\arrowtype<-3
\putdashvector(#1,#2)(#3,#4)#5\arrowtype
\else
\xpos=#1
\ypos=#2
\run=#3
\rise=#4
\arrowlength=#5
\ifnum \arrowtype<0
    \ifnum \run=0
        \advance \ypos by-\arrowlength
    \else
        \tempcounta \arrowlength
        \multiply \tempcounta by\rise
        \divide \tempcounta by\run
        \ifnum\run>0
            \advance \xpos by\arrowlength
            \advance \ypos by\tempcounta
        \else
            \advance \xpos by-\arrowlength
            \advance \ypos by-\tempcounta
        \fi
    \fi
    \multiply \arrowtype by-1
    \multiply \rise by-1
    \multiply \run by-1
\fi
\ifcase \arrowtype
\or \put(\xpos,\ypos){\vector(\run,\rise){\arrowlength}}%
\or \put(\xpos,\ypos){\mvector(\run,\rise)\arrowlength}%
\or \put(\xpos,\ypos){\evector(\run,\rise){\arrowlength}}%
\fi\fi\fi
}}
 
\def\putsplitvector(#1,#2)#3#4{
\xpos #1
\ypos #2
\arrowtype #4
\halflength #3
\arrowlength #3
\gap 140
\advance \halflength by-\gap
\divide \halflength by2
\ifnum\arrowtype>0
   \ifcase \arrowtype
   \or \put(\xpos,\ypos){\line(0,-1){\halflength}}%
       \advance\ypos by-\halflength
       \advance\ypos by-\gap
       \put(\xpos,\ypos){\vector(0,-1){\halflength}}%
   \or \put(\xpos,\ypos){\line(0,-1)\halflength}%
       \put(\xpos,\ypos){\vector(0,-1)3}%
       \advance\ypos by-\halflength
       \advance\ypos by-\gap
       \put(\xpos,\ypos){\vector(0,-1){\halflength}}%
   \or \put(\xpos,\ypos){\line(0,-1)\halflength}%
       \advance\ypos by-\halflength
       \advance\ypos by-\gap
       \put(\xpos,\ypos){\evector(0,-1){\halflength}}%
   \fi
\else \arrowtype=-\arrowtype
   \ifcase\arrowtype
   \or \advance \ypos by-\arrowlength
       \put(\xpos,\ypos){\line(0,1){\halflength}}%
       \advance\ypos by\halflength
       \advance\ypos by\gap
       \put(\xpos,\ypos){\vector(0,1){\halflength}}%
   \or \advance \ypos by-\arrowlength
       \put(\xpos,\ypos){\line(0,1)\halflength}%
       \put(\xpos,\ypos){\vector(0,1)3}%
       \advance\ypos by\halflength
       \advance\ypos by\gap
       \put(\xpos,\ypos){\vector(0,1){\halflength}}%
   \or \advance \ypos by-\arrowlength
       \put(\xpos,\ypos){\line(0,1)\halflength}%
       \advance\ypos by\halflength
       \advance\ypos by\gap
       \put(\xpos,\ypos){\evector(0,1){\halflength}}%
   \fi
\fi
}
 
\def\putmorphism(#1)(#2,#3)[#4`#5`#6]#7#8#9{{%
\run #2
\rise #3
\ifnum\rise=0
  \puthmorphism(#1)[#4`#5`#6]{#7}{#8}#9%
\else\ifnum\run=0
  \putvmorphism(#1)[#4`#5`#6]{#7}{#8}#9%
\else
\setpos(#1)%
\arrowlength #7
\arrowtype #8
\ifnum\run=0
\else\ifnum\rise=0
\else
\ifnum\run>0
    \coefa=1
\else
   \coefa=-1
\fi
\ifnum\arrowtype>0
   \coefb=0
   \coefc=-1
\else
   \coefb=\coefa
   \coefc=1
   \arrowtype=-\arrowtype
\fi
\width=2
\multiply \width by\run
\divide \width by\rise
\ifnum \width<0  \width=-\width\fi
\advance\width by60
\if l#9 \width=-\width\fi
\putbox(\xpos,\ypos){#4}
{\multiply \coefa by\arrowlength
\advance\xpos by\coefa
\multiply \coefa by\rise
\divide \coefa by\run
\advance \ypos by\coefa
\putbox(\xpos,\ypos){#5} }%
{\multiply \coefa by\arrowlength
\divide \coefa by2
\advance \xpos by\coefa
\advance \xpos by\width
\multiply \coefa by\rise
\divide \coefa by\run
\advance \ypos by\coefa
\if l#9%
   \putrbox(\xpos,\ypos){#6}%
\else\if r#9%
   \putlbox(\xpos,\ypos){#6}%
\fi\fi }%
{\multiply \rise by-\coefc
\multiply \run by-\coefc
\multiply \coefb by\arrowlength
\advance \xpos by\coefb
\multiply \coefb by\rise
\divide \coefb by\run
\advance \ypos by\coefb
\multiply \coefc by70
\advance \ypos by\coefc
\multiply \coefc by\run
\divide \coefc by\rise
\advance \xpos by\coefc
\multiply \coefa by140
\multiply \coefa by\run
\divide \coefa by\rise
\advance \arrowlength by\coefa
\ifcase\arrowtype
\or \put(\xpos,\ypos){\vector(\run,\rise){\arrowlength}}%
\or \put(\xpos,\ypos){\mvector(\run,\rise){\arrowlength}}%
\or \put(\xpos,\ypos){\evector(\run,\rise){\arrowlength}}%
\fi}\fi\fi\fi\fi}}

\newcount\numbdashes \newcount\lengthdash \newcount\increment
 
\def\howmanydashes{
\numbdashes=\arrowlength \lengthdash=40
\divide\numbdashes by \lengthdash
\lengthdash=\arrowlength
\divide\lengthdash by \numbdashes
\increment=\lengthdash
\multiply\lengthdash by 3
\divide\lengthdash by 5
}
 
\def\putdashvector(#1)(#2,#3)#4#5{%
\ifnum#3=0 \putdashhvector(#1){#4}#5
\else
\ifnum#2=0
\putdashvvector(#1){#4}#5\fi\fi}
 
\def\putdashhvector(#1,#2)#3#4{{%
\arrowlength=#3 \howmanydashes
\multiput(#1,#2)(\increment,0){\numbdashes}%
{\vrule height .4pt width \lengthdash\unitlength}
\arrowtype=#4 \xpos=#1
\ifnum\arrowtype<0 \advance\arrowtype by 7 \fi
\ifcase\arrowtype
\or \advance\xpos by 10
    \put(\xpos,#2){\vector(-1,0){\lengthdash}}
    \advance\xpos by 40
    \put(\xpos,#2){\vector(-1,0){\lengthdash}}
\or \advance \xpos by 10
    \put(\xpos,#2){\vector(-1,0){\lengthdash}}
    \advance\xpos by  \arrowlength
    \advance\xpos by  -50
    \put(\xpos,#2){\vector(-1,0){\lengthdash}}
\or \advance\xpos by 10
    \put(\xpos,#2){\vector(-1,0){\lengthdash}}
\or \advance\xpos by \arrowlength
    \advance\xpos by -\lengthdash
    \put(\xpos,#2){\vector(1,0){\lengthdash}}
\or {\advance\xpos by 10
    \put(\xpos,#2){\vector(1,0){\lengthdash}}}
    \advance\xpos by \arrowlength
    \advance\xpos by -\lengthdash
    \put(\xpos,#2){\vector(1,0){\lengthdash}}
\or \advance\xpos by \arrowlength
    \advance\xpos by -\lengthdash
    \put(\xpos,#2){\vector(1,0){\lengthdash}}
    \advance\xpos by -40
    \put(\xpos,#2){\vector(1,0){\lengthdash}}
   \fi
}}
 
\def\putdashvvector(#1,#2)#3#4{{%
\arrowlength=#3 \howmanydashes
\ypos=#2 \advance\ypos by -\arrowlength
\multiput(#1,#2)(0,\increment){\numbdashes}%
    {\vrule width .4pt height \lengthdash\unitlength}
\arrowtype=#4 \ypos=#2
\ifnum\arrowtype<0 \advance\arrowtype by 7 \fi
\ifcase\arrowtype
\or \advance\ypos by \arrowlength \advance\ypos by -40
    \put(#1,\ypos){\vector(0,1){\lengthdash}}
    \advance\ypos by -40
    \put(#1,\ypos){\vector(0,1){\lengthdash}}
\or \advance\ypos by 10
    \put(#1,\ypos){\vector(0,1){\lengthdash}}
    \advance\ypos by \arrowlength \advance\ypos by -40
    \put(#1,\ypos){\vector(0,1){\lengthdash}}
\or \advance\ypos by \arrowlength \advance\ypos by -40
    \put(#1,\ypos){\vector(0,1){\lengthdash}}
\or \advance\ypos by 10
    \put(#1,\ypos){\vector(0,-1){\lengthdash}}
\or \advance\ypos by 10
    \put(#1,\ypos){\vector(0,-1){\lengthdash}}
    \advance\ypos by \arrowlength \advance\ypos by -40
    \put(#1,\ypos){\vector(0,-1){\lengthdash}}
\or \advance\ypos by 10
    \put(#1,\ypos){\vector(0,-1){\lengthdash}}
    \advance\ypos by 40
    \put(#1,\ypos){\vector(0,-1){\lengthdash}}
\fi
}}
 
\def\puthmorphism(#1,#2)[#3`#4`#5]#6#7#8{{%
\xpos #1
\ypos #2
\width #6
\arrowlength #6
\arrowtype=#7
\putbox(\xpos,\ypos){#3\vphantom{#4}}%
{\advance \xpos by\arrowlength
\putbox(\xpos,\ypos){\vphantom{#3}#4}}%
\horsize{\tempcounta}{#3}%
\horsize{\tempcountb}{#4}%
\divide \tempcounta by2
\divide \tempcountb by2
\advance \tempcounta by30
\advance \tempcountb by30
\advance \xpos by\tempcounta
\advance \arrowlength by-\tempcounta
\advance \arrowlength by-\tempcountb
\putvector(\xpos,\ypos)(1,0)\arrowlength\arrowtype
\divide \arrowlength by2
\advance \xpos by\arrowlength
\vertsize{\tempcounta}{#5}%
\divide\tempcounta by2
\advance \tempcounta by20
\if a#8 %
   \advance \ypos by\tempcounta
   \putbox(\xpos,\ypos){#5}%
\else
   \advance \ypos by-\tempcounta
   \putbox(\xpos,\ypos){#5}%
\fi}}
 
\def\putvmorphism(#1,#2)[#3`#4`#5]#6#7#8{{%
\xpos #1
\ypos #2
\arrowlength #6
\arrowtype #7
\settowidth{\xlen}{$#5$}%
\putbox(\xpos,\ypos){#3}%
{\advance \ypos by-\arrowlength
\putbox(\xpos,\ypos){#4}}%
{\advance\arrowlength by-140
\advance \ypos by-70
\ifdim\xlen>0pt
   \if m#8%
      \putsplitvector(\xpos,\ypos)\arrowlength\arrowtype
   \else
   \putvector(\xpos,\ypos)(0,-1)\arrowlength\arrowtype
   \fi
\else
   \putvector(\xpos,\ypos)(0,-1)\arrowlength\arrowtype
\fi}%
\ifdim\xlen>0pt
   \divide \arrowlength by2
   \advance\ypos by-\arrowlength
   \if l#8%
      \advance \xpos by-40
      \putrbox(\xpos,\ypos){#5}%
   \else\if r#8%
      \advance \xpos by40
      \putlbox(\xpos,\ypos){#5}%
   \else
      \putbox(\xpos,\ypos){#5}%
   \fi\fi
\fi
}}
 
\def\putsquarep<#1>(#2)[#3;#4`#5`#6`#7]{{%
\setsqparms[#1]%
\setpos(#2)%
\settokens[#3]%
\puthmorphism(\xpos,\ypos)[\tokenc`\tokend`{#7}]{\width}{\arrowtyped}b%
\advance\ypos by \height
\puthmorphism(\xpos,\ypos)[\tokena`\tokenb`{#4}]{\width}{\arrowtypea}a%
\putvmorphism(\xpos,\ypos)[``{#5}]{\height}{\arrowtypeb}l%
\advance\xpos by \width
\putvmorphism(\xpos,\ypos)[``{#6}]{\height}{\arrowtypec}r%
}}
 
\def\putsquare{\@ifnextchar <{\putsquarep}{\putsquarep%
   <\arrowtypea`\arrowtypeb`\arrowtypec`\arrowtyped;\width`\height>}}
\def\square{\@ifnextchar< {\squarep}{\squarep
   <\arrowtypea`\arrowtypeb`\arrowtypec`\arrowtyped;\width`\height>}}
\def\squarep<#1>[#2`#3`#4`#5;#6`#7`#8`#9]{{
\setsqparms[#1]
\diagram
\putsquarep<\arrowtypea`\arrowtypeb`\arrowtypec`
\arrowtyped;\width`\height>
(0,0)[#2`#3`#4`{#5};#6`#7`#8`{#9}]
\enddiagram
}}                                                 
\def\putptrianglep<#1>(#2,#3)[#4`#5`#6;#7`#8`#9]{{%
\settriparms[#1]%
\xpos=#2 \ypos=#3
\advance\ypos by \height
\puthmorphism(\xpos,\ypos)[#4`#5`{#7}]{\height}{\arrowtypea}a%
\putvmorphism(\xpos,\ypos)[`#6`{#8}]{\height}{\arrowtypeb}l%
\advance\xpos by\height
\putmorphism(\xpos,\ypos)(-1,-1)[``{#9}]{\height}{\arrowtypec}r%
}}
 
\def\putptriangle{\@ifnextchar <{\putptrianglep}{\putptrianglep
   <\arrowtypea`\arrowtypeb`\arrowtypec;\height>}}
\def\ptriangle{\@ifnextchar <{\ptrianglep}{\ptrianglep
   <\arrowtypea`\arrowtypeb`\arrowtypec;\height>}}
\def\ptrianglep<#1>[#2`#3`#4;#5`#6`#7]{{
\settriparms[#1]
\diagram
\putptrianglep<\arrowtypea`\arrowtypeb`
\arrowtypec;\height>
(0,0)[#2`#3`#4;#5`#6`{#7}]
\enddiagram
}}                                            
 
\def\putqtrianglep<#1>(#2,#3)[#4`#5`#6;#7`#8`#9]{{%
\settriparms[#1]%
\xpos=#2 \ypos=#3
\advance\ypos by\height
\puthmorphism(\xpos,\ypos)[#4`#5`{#7}]{\height}{\arrowtypea}a%
\putmorphism(\xpos,\ypos)(1,-1)[``{#8}]{\height}{\arrowtypeb}l%
\advance\xpos by\height
\putvmorphism(\xpos,\ypos)[`#6`{#9}]{\height}{\arrowtypec}r%
}}
 
\def\putqtriangle{\@ifnextchar <{\putqtrianglep}{\putqtrianglep
   <\arrowtypea`\arrowtypeb`\arrowtypec;\height>}}
\def\qtriangle{\@ifnextchar <{\qtrianglep}{\qtrianglep
   <\arrowtypea`\arrowtypeb`\arrowtypec;\height>}}
\def\qtrianglep<#1>[#2`#3`#4;#5`#6`#7]{{
\settriparms[#1]
\width=\height                                
\diagram
\putqtrianglep<\arrowtypea`\arrowtypeb`
\arrowtypec;\height>
(0,0)[#2`#3`#4;#5`#6`{#7}]
\enddiagram
}}
 
\def\putdtrianglep<#1>(#2,#3)[#4`#5`#6;#7`#8`#9]{{%
\settriparms[#1]%
\xpos=#2 \ypos=#3
\puthmorphism(\xpos,\ypos)[#5`#6`{#9}]{\height}{\arrowtypec}b%
\advance\xpos by \height \advance\ypos by\height
\putmorphism(\xpos,\ypos)(-1,-1)[``{#7}]{\height}{\arrowtypea}l%
\putvmorphism(\xpos,\ypos)[#4``{#8}]{\height}{\arrowtypeb}r%
}}
 
\def\putdtriangle{\@ifnextchar <{\putdtrianglep}{\putdtrianglep
   <\arrowtypea`\arrowtypeb`\arrowtypec;\height>}}
\def\dtriangle{\@ifnextchar <{\dtrianglep}{\dtrianglep
   <\arrowtypea`\arrowtypeb`\arrowtypec;\height>}}
\def\dtrianglep<#1>[#2`#3`#4;#5`#6`#7]{{
\settriparms[#1]
\width=\height                                
\diagram
\putdtrianglep<\arrowtypea`\arrowtypeb`
\arrowtypec;\height>
(0,0)[#2`#3`#4;#5`#6`{#7}]
\enddiagram
}}
 
\def\putbtrianglep<#1>(#2,#3)[#4`#5`#6;#7`#8`#9]{{%
\settriparms[#1]%
\xpos=#2 \ypos=#3
\puthmorphism(\xpos,\ypos)[#5`#6`{#9}]{\height}{\arrowtypec}b%
\advance\ypos by\height
\putmorphism(\xpos,\ypos)(1,-1)[``{#8}]{\height}{\arrowtypeb}r%
\putvmorphism(\xpos,\ypos)[#4``{#7}]{\height}{\arrowtypea}l%
}}
 
\def\putbtriangle{\@ifnextchar <{\putbtrianglep}{\putbtrianglep
   <\arrowtypea`\arrowtypeb`\arrowtypec;\height>}}
\def\btriangle{\@ifnextchar <{\btrianglep}{\btrianglep
   <\arrowtypea`\arrowtypeb`\arrowtypec;\height>}}
\def\btrianglep<#1>[#2`#3`#4;#5`#6`#7]{{
\settriparms[#1]
\width=\height                               
\diagram
\putbtrianglep<\arrowtypea`\arrowtypeb`
\arrowtypec;\height>
(0,0)[#2`#3`#4;#5`#6`{#7}]
\enddiagram
}}
 
\def\putAtrianglep<#1>(#2,#3)[#4`#5`#6;#7`#8`#9]{{%
\settriparms[#1]%
\xpos=#2 \ypos=#3
{\multiply \height by2
\puthmorphism(\xpos,\ypos)[#5`#6`{#9}]{\height}{\arrowtypec}b}%
\advance\xpos by\height \advance\ypos by\height
\putmorphism(\xpos,\ypos)(-1,-1)[#4``{#7}]{\height}{\arrowtypea}l%
\putmorphism(\xpos,\ypos)(1,-1)[``{#8}]{\height}{\arrowtypeb}r%
}}
 
\def\putAtriangle{\@ifnextchar <{\putAtrianglep}{\putAtrianglep
   <\arrowtypea`\arrowtypeb`\arrowtypec;\height>}}
\def\Atriangle{\@ifnextchar <{\Atrianglep}{\Atrianglep
   <\arrowtypea`\arrowtypeb`\arrowtypec;\height>}}
\def\Atrianglep<#1>[#2`#3`#4;#5`#6`#7]{{
\settriparms[#1]
\width=\height                                     
\diagram
\putAtrianglep<\arrowtypea`\arrowtypeb`
\arrowtypec;\height>
(0,0)[#2`#3`#4;#5`#6`{#7}]
\enddiagram
}}
 
\def\putAtrianglepairp<#1>(#2)[#3;#4`#5`#6`#7`#8]{{%
\settripairparms[#1]%
\setpos(#2)%
\settokens[#3]%
\puthmorphism(\xpos,\ypos)[\tokenb`\tokenc`{#7}]{\height}{\arrowtyped}b%
\advance\xpos by\height
\puthmorphism(\xpos,\ypos)[\phantom{\tokenc}`\tokend`{#8}]%
{\height}{\arrowtypee}b%
\advance\ypos by\height
\putmorphism(\xpos,\ypos)(-1,-1)[\tokena``{#4}]{\height}{\arrowtypea}l%
\putvmorphism(\xpos,\ypos)[``{#5}]{\height}{\arrowtypeb}m%
\putmorphism(\xpos,\ypos)(1,-1)[``{#6}]{\height}{\arrowtypec}r%
}}
 
\def\putAtrianglepair{\@ifnextchar <{\putAtrianglepairp}{\putAtrianglepairp%
   <\arrowtypea`\arrowtypeb`\arrowtypec`\arrowtyped`\arrowtypee;\height>}}
\def\Atrianglepair{\@ifnextchar <{\Atrianglepairp}{\Atrianglepairp%
   <\arrowtypea`\arrowtypeb`\arrowtypec`\arrowtyped`\arrowtypee;\height>}}
 
\def\Atrianglepairp<#1>[#2;#3`#4`#5`#6`#7]{{
\settripairparms[#1]
\settokens[#2]
\width=\height                                
\diagram
\putAtrianglepairp                            
<\arrowtypea`\arrowtypeb`\arrowtypec`
\arrowtyped`\arrowtypee;\height>
(0,0)[{#2};#3`#4`#5`#6`{#7}]
\enddiagram
}}
 
\def\putVtrianglep<#1>(#2,#3)[#4`#5`#6;#7`#8`#9]{{%
\settriparms[#1]%
\xpos=#2 \ypos=#3
\advance\ypos by\height
{\multiply\height by2
\puthmorphism(\xpos,\ypos)[#4`#5`{#7}]{\height}{\arrowtypea}a}%
\putmorphism(\xpos,\ypos)(1,-1)[`#6`{#8}]{\height}{\arrowtypeb}l%
\advance\xpos by\height
\advance\xpos by\height
\putmorphism(\xpos,\ypos)(-1,-1)[``{#9}]{\height}{\arrowtypec}r%
}}
 
\def\putVtriangle{\@ifnextchar <{\putVtrianglep}{\putVtrianglep
   <\arrowtypea`\arrowtypeb`\arrowtypec;\height>}}
\def\Vtriangle{\@ifnextchar <{\Vtrianglep}{\Vtrianglep
   <\arrowtypea`\arrowtypeb`\arrowtypec;\height>}}
\def\Vtrianglep<#1>[#2`#3`#4;#5`#6`#7]{{
\settriparms[#1]
\width=\height                                 
\diagram
\putVtrianglep<\arrowtypea`\arrowtypeb`
\arrowtypec;\height>
(0,0)[#2`#3`#4;#5`#6`{#7}]
\enddiagram
}}
 
\def\putVtrianglepairp<#1>(#2)[#3;#4`#5`#6`#7`#8]{{
\settripairparms[#1]%
\setpos(#2)%
\settokens[#3]%
\advance\ypos by\height
\putmorphism(\xpos,\ypos)(1,-1)[`\tokend`{#6}]{\height}{\arrowtypec}l%
\puthmorphism(\xpos,\ypos)[\tokena`\tokenb`{#4}]{\height}{\arrowtypea}a%
\advance\xpos by\height
\puthmorphism(\xpos,\ypos)[\phantom{\tokenb}`\tokenc`{#5}]%
{\height}{\arrowtypeb}a%
\putvmorphism(\xpos,\ypos)[``{#7}]{\height}{\arrowtyped}m%
\advance\xpos by\height
\putmorphism(\xpos,\ypos)(-1,-1)[``{#8}]{\height}{\arrowtypee}r%
}}
 
\def\putVtrianglepair{\@ifnextchar <{\putVtrianglepairp}{\putVtrianglepairp%
    <\arrowtypea`\arrowtypeb`\arrowtypec`\arrowtyped`\arrowtypee;\height>}}
\def\Vtrianglepair{\@ifnextchar <{\Vtrianglepairp}{\Vtrianglepairp%
    <\arrowtypea`\arrowtypeb`\arrowtypec`\arrowtyped`\arrowtypee;\height>}}
\def\Vtrianglepairp<#1>[#2;#3`#4`#5`#6`#7]{{
\settripairparms[#1]
\settokens[#2]
\diagram
\putVtrianglepairp                             
<\arrowtypea`\arrowtypeb`\arrowtypec`
\arrowtyped`\arrowtypee;\height>
(0,0)[{#2};#3`#4`#5`#6`{#7}]
\enddiagram
}}

\def\putCtrianglep<#1>(#2,#3)[#4`#5`#6;#7`#8`#9]{{%
\settriparms[#1]%
\xpos=#2 \ypos=#3
\advance\ypos by\height
\putmorphism(\xpos,\ypos)(1,-1)[``{#9}]{\height}{\arrowtypec}l%
\advance\xpos by\height
\advance\ypos by\height
\putmorphism(\xpos,\ypos)(-1,-1)[#4`#5`{#7}]{\height}{\arrowtypea}l%
{\multiply\height by 2
\putvmorphism(\xpos,\ypos)[`#6`{#8}]{\height}{\arrowtypeb}r}%
}}
 
\def\putCtriangle{\@ifnextchar <{\putCtrianglep}{\putCtrianglep
    <\arrowtypea`\arrowtypeb`\arrowtypec;\height>}}
\def\Ctriangle{\@ifnextchar <{\Ctrianglep}{\Ctrianglep
    <\arrowtypea`\arrowtypeb`\arrowtypec;\height>}}
\def\Ctrianglep<#1>[#2`#3`#4;#5`#6`#7]{{
\settriparms[#1]
\width=\height                               
\diagram
\putCtrianglep<\arrowtypea`\arrowtypeb`
\arrowtypec;\height>
(0,0)[#2`#3`#4;#5`#6`{#7}]
\enddiagram
}}                                           
\def\putDtrianglep<#1>(#2,#3)[#4`#5`#6;#7`#8`#9]{{%
\settriparms[#1]%
\xpos=#2 \ypos=#3
\advance\xpos by\height \advance\ypos by\height
\putmorphism(\xpos,\ypos)(-1,-1)[``{#9}]{\height}{\arrowtypec}r%
\advance\xpos by-\height \advance\ypos by\height
\putmorphism(\xpos,\ypos)(1,-1)[`#5`{#8}]{\height}{\arrowtypeb}r%
{\multiply\height by 2
\putvmorphism(\xpos,\ypos)[#4`#6`{#7}]{\height}{\arrowtypea}l}%
}}
 
\def\putDtriangle{\@ifnextchar <{\putDtrianglep}{\putDtrianglep
    <\arrowtypea`\arrowtypeb`\arrowtypec;\height>}}
\def\Dtriangle{\@ifnextchar <{\Dtrianglep}{\Dtrianglep
   <\arrowtypea`\arrowtypeb`\arrowtypec;\height>}}
\def\Dtrianglep<#1>[#2`#3`#4;#5`#6`#7]{{
\settriparms[#1]
\width=\height                              
\diagram
\putDtrianglep<\arrowtypea`\arrowtypeb`
\arrowtypec;\height>
(0,0)[#2`#3`#4;#5`#6`{#7}]
\enddiagram
}}                                          
\def\setrecparms[#1`#2]{\width=#1 \height=#2}%
 
\def\recursep<#1`#2>[#3;#4`#5`#6`#7`#8]{{%
\width=#1 \height=#2
\settokens[#3]
\settowidth{\tempdimen}{$\tokena$}
\ifdim\tempdimen=0pt
  \savebox{\tempboxa}{\hbox{$\tokenb$}}%
  \savebox{\tempboxb}{\hbox{$\tokend$}}%
  \savebox{\tempboxc}{\hbox{$#6$}}%
\else
  \savebox{\tempboxa}{\hbox{$\hbox{$\tokena$}\times\hbox{$\tokenb$}$}}%
  \savebox{\tempboxb}{\hbox{$\hbox{$\tokena$}\times\hbox{$\tokend$}$}}%
  \savebox{\tempboxc}{\hbox{$\hbox{$\tokena$}\times\hbox{$#6$}$}}%
\fi
\ypos=\height
\divide\ypos by 2
\xpos=\ypos
\advance\xpos by \width
\bfig
\putCtrianglep<-1`1`1;\ypos>(0,0)[`\tokenc`;#5`#6`{#7}]%
\puthmorphism(\ypos,0)[\tokend`\usebox{\tempboxb}`{#8}]{\width}{-1}b%
\puthmorphism(\ypos,\height)[\tokenb`\usebox{\tempboxa}`{#4}]{\width}{-1}a%
\advance\ypos by \width
\putvmorphism(\ypos,\height)[``\usebox{\tempboxc}]{\height}1r%
\efig
}}
 
\def\recurse{\@ifnextchar <{\recursep}{\recursep<\width`\height>}}
 
\def\puttwohmorphisms(#1,#2)[#3`#4;#5`#6]#7#8#9{{%
\puthmorphism(#1,#2)[#3`#4`]{#7}0a
\ypos=#2
\advance\ypos by 20
\puthmorphism(#1,\ypos)[\phantom{#3}`\phantom{#4}`#5]{#7}{#8}a
\advance\ypos by -40
\puthmorphism(#1,\ypos)[\phantom{#3}`\phantom{#4}`#6]{#7}{#9}b
}}
 
\def\puttwovmorphisms(#1,#2)[#3`#4;#5`#6]#7#8#9{{%
\putvmorphism(#1,#2)[#3`#4`]{#7}0a
\xpos=#1
\advance\xpos by -20
\putvmorphism(\xpos,#2)[\phantom{#3}`\phantom{#4}`#5]{#7}{#8}l
\advance\xpos by 40
\putvmorphism(\xpos,#2)[\phantom{#3}`\phantom{#4}`#6]{#7}{#9}r
}}
 
\def\puthcoequalizer(#1)[#2`#3`#4;#5`#6`#7]#8#9{{%
\setpos(#1)%
\puttwohmorphisms(\xpos,\ypos)[#2`#3;#5`#6]{#8}11%
\advance\xpos by #8
\puthmorphism(\xpos,\ypos)[\phantom{#3}`#4`#7]{#8}1{#9}
}}
 
\def\putvcoequalizer(#1)[#2`#3`#4;#5`#6`#7]#8#9{{%
\setpos(#1)%
\puttwovmorphisms(\xpos,\ypos)[#2`#3;#5`#6]{#8}11%
\advance\ypos by -#8
\putvmorphism(\xpos,\ypos)[\phantom{#3}`#4`#7]{#8}1{#9}
}}
 
\def\putthreehmorphisms(#1)[#2`#3;#4`#5`#6]#7(#8)#9{{%
\setpos(#1) \settypes(#8)
\if a#9 %
     \vertsize{\tempcounta}{#5}%
     \vertsize{\tempcountb}{#6}%
     \ifnum \tempcounta<\tempcountb \tempcounta=\tempcountb \fi
\else
     \vertsize{\tempcounta}{#4}%
     \vertsize{\tempcountb}{#5}%
     \ifnum \tempcounta<\tempcountb \tempcounta=\tempcountb \fi
\fi
\advance \tempcounta by 60
\puthmorphism(\xpos,\ypos)[#2`#3`#5]{#7}{\arrowtypeb}{#9}
\advance\ypos by \tempcounta
\puthmorphism(\xpos,\ypos)[\phantom{#2}`\phantom{#3}`#4]{#7}{\arrowtypea}{#9}
\advance\ypos by -\tempcounta \advance\ypos by -\tempcounta
\puthmorphism(\xpos,\ypos)[\phantom{#2}`\phantom{#3}`#6]{#7}{\arrowtypec}{#9}
}}
 
\def\setarrowtoks[#1`#2`#3`#4`#5`#6]{%
\def\toka{#1}
\def\tokb{#2}
\def\tokc{#3}
\def\tokd{#4}
\def\toke{#5}
\def\tokf{#6}
}
\def\hex{\@ifnextchar <{\hexp}{\hexp<1000`400>}}
\def\hexp<#1`#2>[#3`#4`#5`#6`#7`#8;#9]{%
\setarrowtoks[#9]
\yext=#2 \advance \yext by #2
\xext=#1 \advance\xext by \yext
\bfig
\putCtriangle<-1`0`1;#2>(0,0)[`#5`;\tokb``\tokd]
\xext=#1 \yext=#2 \advance \yext by #2
\putsquare<1`0`0`1;\xext`\yext>(#2,0)[#3`#4`#7`#8;\toka```\tokf]
\advance \xext by #2
\putDtriangle<0`1`-1;#2>(\xext,0)[`#6`;`\tokc`\toke]
\efig
}

\usepackage{amscd}
\usepackage{amssymb}

\newtheorem{Theorem}{Theorem}[section]
\newtheorem{Lemma}[Theorem]{Lemma}
\newtheorem{Proposition}[Theorem]{Proposition}
\newtheorem{Definition}[Theorem]{Definition}
\newtheorem{Corollary}[Theorem]{Corollary}

\newenvironment{Remark}
  {\begin{flushleft}\textbf{Remark.}\begin{sl} }
  {\end{sl}\end{flushleft}}

\def\ot{\otimes}
\def\ep{\varepsilon}

\def\Hom{\operatorname{Hom}}
\def\Reg{\operatorname{Reg}}

\def\Aut{\operatorname{Aut}}
\def\Bimeas{\operatorname{Bimeas}}
\def\Hopf{\operatorname{Hopf}}
\def\Bialg{\operatorname{Bialg}}
\def\Alg{\operatorname{Alg}}
\def\Vect{\operatorname{Vect}}
\def\Coalg{\operatorname{Coalg}}

\def\m{\operatorname{m}}

\def\B{B}

\def\io{\iota}

\def\De{\Delta}

\def\om{\omega}

\def\id{\operatorname{id}}
\def\ol{\overline}

\title{On Bimeasurings}
\author{L. Grunenfelder and M. Mastnak}
\address{Department of Mathematics and Statistics,
Dalhousie University, Halifax, Nova Scotia, Canada, B3H 3J5}
\email{luzius@mathstat.dal.ca,\ mastnak@mathstat.dal.ca}
\thanks{}
\date{}

\begin{document}

\tolerance=1000

\begin{abstract}
We introduce and study bimeasurings from pairs of bialgebras to algebras.
It is shown that the universal bimeasuring bialgebra construction, which arises
from Sweedler's universal measuring coalgebra construction and
generalizes the finite dual, gives rise to a contravariant
functor on the category of bialgebras adjoint to itself. An interpretation of
bimeasurings as algebras in the category of Hopf modules is considered.
\end{abstract}

\maketitle

\setcounter{section}{-1}

\section{Introduction}

Measurings have first been introduced and studied my M.E. Sweedler
\cite{Sw}. They correspond to homomorphisms of algebras over a
coalgebra which are cofree as comodules \cite{GP}. There is a
universal measuring coalgebra $M(B,A)$ and measuring $\theta\colon
M(B,A)\ot B\to A$ for every pair of algebras $A$ and $B$ such that
$C$-measurings from $B$ to $A$ correspond bijectively to coalgebra
maps from $C$ to $M(B,A)$. If $B$ is a Hopf algebra and $A$ is
commutative then $M(B,A)$ carries a natural Hopf algebra structure
\cite{Ma}. If in addition, $C$ is a Hopf algebra then one may
consider maps $\psi\colon C\ot B\to A$ which measure in both
variables $C$ and $B$. In the cocommutative case these
bimeasurings account for the {\lq\lq mixed\rq\rq} term in the
second Sweedler cohomology group
$$
H^2(C\ot B,A)\simeq H^2(C,A)\oplus H^2(B,A)\oplus P(B,C,A)
$$
as shown in \cite{Ma}. If $A$ is commutative then universal bimeasuring Hopf
algebras (and universal bimeasuring) $B(C,A)$ and $B(B,A)$ exist so that
bimeasurings $\theta\colon C\ot B\to A$ bijectively
correspond to Hopf algebra maps from $C$ to $B(B,A)$ as well as Hopf algebra
maps from $B$ to $B(C,A)$. In fact
$$
\Hopf(C,B(B,A))\simeq \Bimeas(C\ot B,A)\simeq \Hopf(B,B(C,A))
$$
and hence the functor $B(\_,A)$ on the category of Hopf algebras
is adjoint to itself. In the special case $A=k$ this gives a new
proof that the finite dual construction $\_^\circ=B(\_,k)$ is
adjoint to itself \cite{Ta} (see \cite{Mi} for the proof).
Moreover, there is a natural injective map
$$
B(C,A)\ot B(B,A)\to B(C\ot B,A),
$$
which is always an isomorphism in the cocommutative case, and restricts to the
well known isomorphism $C^\circ\ot B^\circ\simeq (C\ot B)^\circ$ when $A=k$.

There is a natural notion of bimeasuring from an abelian matched
pair of Hopf algebras $H=C\bowtie B$ to a commutative algebra $A$
extending that of ordinary cocommutative bimeasurings. These
skew-bimeasurings form an abelian group under convolution
isomorphic to the first matched pair cohomology group
$\mathcal{H}^1(C,B,A)$ with coefficients in $A$ described in
\cite{GM}. This group also corresponds to a subgroup of the group
of $A$-linear automorphisms of the trivial $H$-comodule $H\ot A$
and thus to a group of Hopf algebra structures on $H\ot A$, each
making $H\ot A$ an algebra in the category of Hopf modules.

\section{Preliminaries}
\subsection{Notation}
All vector spaces (algebras, coalgebras, bialgebras)
will be over a ground field $k$. If $A$ is an algebra and $C$ a coalgebra, then
$\Hom(C,A)$ denotes the convolution algebra of all linear maps from $C$ to $A$.
The unit and the multiplication on $A$ are denoted by $\eta\colon k\to A$ and $\m\colon A\ot A\to A$;
the counit and the comultiplication on $C$ are denoted by $\ep\colon C\to k$ and $\De\colon C\to C\ot C$.
We use Sweedler's sigma notation for comultiplication: $\De(c)=c_1\ot c_2$,
$(1\ot\De)\De(c)=c_1\ot c_2\ot c_3$ etc. If $f\colon U\ot V\to W$ is a linear map than we often write
$f(u,v)$ instead of $f(u\ot v)$.

\subsection{Abelianization}
Let $H$ be an algebra and $I\subseteq H$ the algebra ideal generated by all
commutators, i.e. all elements of the form $[x,y]=xy-yx$. If $H$ is a Hopf algebra
(bialgebra) then $I$ is a Hopf ideal (biideal). This is easily observed by the following
identities:
\begin{eqnarray*}
S[x,y] &=& [S(y),S(x)],\\
\De [x,y] &=& x_1y_1\ot x_2y_2 - y_1x_1\ot y_2x_2\\
&=& [x_1,y_1]\ot x_2y_2 + y_1x_1\ot[x_2,y_2].
\end{eqnarray*}
We call the quotient algebra (Hopf algebra, bialgebra)
$H_{ab}=H/I$ the abelianization of $H$. It is the largest
commutative quotient of $H$ in the sense that if $K$ is a
commutative algebra (bialgebra) and $f\colon H\to K$ is an algebra
(bialgebra) map, then there exists a unique algebra (bialgebra)
map $\overline{f}\colon H_{ab}\to K$ such that
$f=\overline{f}\pi$, where $\pi\colon H\to H_{ab}$ is the
canonical projection.

If $H$ is a Hopf algebra then $I$ is also the Hopf ideal generated by
$\langle [x,y]_H-\ep(xy)\rangle$, where $[x,y]_H=S(x_1)S(y_1)x_2y_2$.

\subsection{Cocommutative part}
For a coalgebra $H$ we define $H_c$, the cocommutative part of H,
to be the largest cocommutative subcoalgebra of $H$ (it is
obtained as a sum of all cocommutative subcoalgebras of $H$, hence
it always exists). If $H$ is a bialgebra, then $H_c$ is a
bialgebra as well (the algebra generated by $H_c$ is also a
cocommutative subcoalgebra of $H$ and must therefore be equal to
$H_c$). Finally, if $H$ is a Hopf algebra, then so is $H_c$. This
is seen by noting that $S(H_c)$ is also a cocommutative
subcoalgebra of $H$. If $f\colon K\to H$ is a coalgebra
(bialgebra) map and $K$ is cocommutative, then clearly
$f(K)\subseteq H_c$; in other words, there exists a unique
coalgebra (bialgebra) map $\overline{f}\colon K\to H_c$ such that
$f=\iota\overline{f}$ (here $\iota\colon H_c\to H$ is the obvious
map).

\subsection{Measuring}
Let $A$, $B$, be algebras, $C$ a coalgebra.
\begin{Proposition}[\cite{Sw}, 7.0.1]
A map $\psi\colon C\ot B\to A$ corresponds to an algebra map $\rho\colon
B\to\Hom(C,A)$, $\rho(b)(c)=\psi(c, b)$ if and only if
\begin{enumerate}
\item $\psi(c, bb')=\psi(c_1, b)\psi(c_2, b'),$ \item
$\psi(c, 1)=\ep(c)$
\end{enumerate}
\end{Proposition}
If the equivalent conditions from the proposition above are
satisfied, we say that $\psi$ is a {\bf measuring}, or that $C$
{\bf measures} $B$ to $A$.

\begin{Theorem}[\cite{Sw}, 7.0.4]
If $A$ and $B$ are algebras then there exits a unique measuring
$\theta\colon M\ot B\to A$ so that for any measuring $f\colon C\ot
B\to A$ there exists a unique coalgebra map $\overline{f}\colon
C\to M$, s.t. $f=\theta(\overline{f}\ot 1)$.
\end{Theorem}
The measuring $\theta\colon M\ot B\to A$ from the theorem above is
called the \textbf{universal measuring} and the coalgebra
$M=M(B,A)$ the \textbf{universal measuring coalgebra}. The functor
$M(\_,A)\colon \Alg^{op}\to \Coalg$ is right adjoint to
$\Hom(\_,A) \colon \Coalg\to \Alg^{op}$. In particular, if $A=k$
then $M(B,A) = M(B,k) = B^\circ$ (the finite dual) and if $B=k$
then $M(B,A)= M(k,A)=k$.

In the construction of the universal bimeasurings, we shall use
the following technical lemma.
\begin{Lemma}\label{l1}
Let $A$ and $B$ be algebras and $\psi\colon M\ot B\to A$ the
universal measuring. If $C$ is a coalgebra and $f$ and $g$
coalgebra maps from $C$ to $M$ such that $\theta(f\ot 1_B)=
\theta(g\ot 1_B)$, then $f=g$.
\end{Lemma}
\begin{proof}
Observe that $\theta(f\ot 1)=\theta(g\ot 1)\colon C\ot B\to A$ is
a measuring and hence by the universal property we have $f=g$.
\end{proof}

If we restrict ourselves to the category of cocommutative
coalgebras, then we talk about \textbf{universal cocommutative
measurings} and \textbf{universal cocommutative measuring
coalgebras}. These were considered in \cite{GP}. In this case, if $C$ is
cocommutative, then $C$-measurings $\psi\colon C\ot B\to A$ are in
bijective correspondence with $C$-algebra maps $\chi\colon C\ot B\to C\ot A$,
given by $\chi=(1\ot\psi)(\De\ot 1)$ and $\psi=(\ep\ot 1)\chi$.

\begin{Proposition}
If $A$ and $B$ are algebras, then the universal cocommutative
measuring coalgebra $M_c(B,A)$ is isomorphic to the cocommutative
part $M(B,A)_c$ of the universal measuring coalgebra $M(B,A)$.
\end{Proposition}
\begin{proof}
Note that $M(B,A)_c$ has the required universal property.
\end{proof}

\section{Bimeasuring}

\begin{Definition}
If $N$ and $T$ are bialgebras and $A$ an algebra, then a map $\psi\colon
N\ot T\to A$ is a \textbf{bimeasuring} if $N$ measures $T$ to $A$ and $T$ measures
$N$ to $A$, i.e.
\begin{eqnarray*}
\psi(nm, t)=\psi(n, t_1)\psi(m, t_2),\; \psi(1_N, t)=\ep(t)\\
\psi(n, ts)=\psi(n_1, t)\psi(n_2, s),\; \psi(n,
1_T)=\ep(n).
\end{eqnarray*}
for $n,m\in N$ and $t,s\in T$.
\end{Definition}

\begin{Definition}
Let $T$ be a bialgebra and $A$ an algebra. If a
bimeasuring $\theta\colon B\ot T\to A$ is such that for every
bimeasuring $f\colon N\ot T\to A$, there exists a unique bialgebra
map $\overline{f}\colon N\to B$ with the property
$f=(\overline{f}\ot 1)\theta$, then $\theta$ is called the (left)
\textbf{universal bimeasuring} and $B=B(T,A)$ is called the (left)
\textbf{universal bimeasuring bialgebra}.

If we limit ourselves to cocommutative $B's$ and $N's$ we talk
about the \textbf{universal cocommutative bimeasurings} and we
denote the \textbf{universal cocommutative bimeasuring bialgebra}
(if it exists) by $B_c(T,A)$.
\end{Definition}

\subsection{Bimeasurings over commutative algebras} The following proposition shows,
that universal bimeasurings exist whenever the algebra $A$ is
commutative.
\begin{Proposition}
If $T$ is a bialgebra, $A$ a commutative algebra, and
$\theta\colon M\ot T\to A$ the universal measuring, then there
exists a unique algebra structure on $M$ so that $T$ measures $M$
to $A$, i.e. $\theta(fg, t)=\theta(f, t_1)\theta(g, t_2)$
and $\theta(1_{M}, t)=\ep(t)$.

Furthermore, with this algebra structure $M$ becomes a bialgebra
and $\theta$ the universal bimeasuring. If $T$ is a Hopf algebra,
then so is $M$.
\end{Proposition}
\begin{proof}
Observe that $\om\colon M\ot M\ot T\to A$, given by $\om(m\ot
m', t)=\theta(m, t_1)\theta(m', t_2)$ is a measuring and
define the multiplication $\m\colon M\ot M\to M$, to be the unique
coalgebra map so that $\theta(\m\ot 1)=\om$.

Similarly the unit $\eta\colon k\to M$ is the unique coalgebra map
so that $\theta(\eta\ot 1)=\eta_A\ep_N$.

The associativity and the unit conditions follow from Lemma
\ref{l1} by noting that $\theta(\m(\m\ot 1_M)\ot
1_T)=\theta(\m(1_M\ot \m)\ot 1_T)$, $\theta(\m(\eta\ot
1_M)\ot 1_T)=\theta(\tau_l\ot 1_T)$ and $\theta(\m(1_M\ot \eta)\ot
1_T)=\theta(\tau_r\ot 1_T)$ (here $\tau_l$ and $\tau_r$ denote the
canonical isomorphisms from $k\ot M$ and $M\ot k$ respectively to
$M$).

Since the multiplication and the unit are coalgebra maps, $M$ must be a
bialgebra and $\theta\colon M\ot T\to A$ a bimeasuring.
We claim $\theta$ is the universal bimeasuring. Let $f\colon N\ot T\to A$ be a
bimeasuring. Since $f$ is a measuring, there exists a unique
coalgebra map $\overline{f}\colon N\to M$ so that
$\psi(\overline{f}\ot 1)=f$. It remains to show that $f$ is also
an algebra map. This follows from Lemma \ref{l1}, since we have
$\psi(\m(\overline{f}\ot \overline{f})\ot
1)=\psi(\overline{f}\m\ot 1)$ and $\psi(\overline{f}\eta_N\ot
1)=\psi(\eta_M\ot 1)$.

We conclude by pointing out that if $T$ is a Hopf algebra, then
the unique coalgebra map $S\colon M\to M^{cop}$ (here $M^{cop}$
denotes the coopposite coalgebra of $M$) satisfying $\theta(S\ot
1)=\theta(1\ot S_N)$, defines the antipode on $M$.
\end{proof}

\begin{Theorem}
If $A$ is a commutative algebra, then the universal bimeasuring
bialgebra construction gives rise to a contravariant functor
$B(\_,A)$ on the category of bialgebras that is adjoint to itself.
The functor restricts to Hopf algebras.
\end{Theorem}
\begin{proof}
It is easy to see that the construction is functorial.

Let $T$ and $N$ be bialgebras. We shall display a canonical
bijection $$\psi_{T,N}\colon \Bialg(T,B(N,A))\to
\Bialg(N,B(T,A)).$$ It is observed from the diagram below.
$$
\bfig \putmorphism(0,750)(3,-2)[{B(T,A)\ot
T}`{A}`{\theta_T}]{1000}{1}l
\puthmorphism(0,750)[\phantom{B(T,A)\ot T}`{N\ot T}`{1\ot
\overline{f}}]{1000}{-1}a \puthmorphism(1000,750)[\phantom{N\ot
T}`{N\ot B(N,A)}`{1\ot f}]{1000}{1}a
\putmorphism(2000,750)(-3,-2)[{N\ot
B(N,A)}`\phantom{A}`{\theta_N}]{1000}{1}r
\putvmorphism(1000,750)[\phantom{N\ot T}`\phantom{A}`{}]{700}{1}l
\efig
$$
More precisely, if $\theta_T\colon B(T,A)\ot T\to A$ and
$\theta_N\colon N\ot B(N,A)\to A$ are universal bimeasurings and
$f\colon T\to B(N,A)$ is a bialgebra map, then $\theta_N(1\ot
f)\colon N\ot T\to A$ is a bimeasuring and we define
$\psi_{T,N}(f)=\overline{f}\colon N\to B(T,A)$ to be the unique
bialgebra map such that $\theta_T(\overline{f}\ot 1)=\theta_N(1\ot
f)\colon N\ot T\to A$. If $g\colon S\to B(T,A)$ is a bialgebra
map, than define $\xi_{N,T}(g)=\overline{g}\colon T\to B(N,A)$ to
be the unique bialgebra map so that
$\theta_T(1\ot\overline{g})=\theta_N(g\ot 1)$ and note that
$\xi_{N,T}$ is the inverse of $\psi_{T,N}$.

We shall conclude the proof by showing that
$\psi_{R,N}(f\alpha)=B(\alpha,A)\psi_{T,N}(f)$, if $\alpha\colon
R\to T$ is a bialgebra map. Indeed, if $\theta_R\colon B(R,A)\ot
R\to A$ is the universal bimeasuring, then
$\theta_R(B(\alpha,A)\overline{f}\ot 1)=\theta_R(B(\alpha,A)\ot
1)(\overline{f}\ot 1)=\theta_T(1\ot \alpha)(\overline{f}\ot 1)=
\theta_T(\overline{f}\ot 1)(1\ot \alpha)=\theta_N(1\ot f)(1\ot
\alpha)=\theta_N(1\ot f\alpha)=\theta_R(\psi_{R,N}(f\alpha)\ot
1)$. Hence we are done by Lemma \ref{l1}.
\end{proof}

\begin{Corollary}\cite{Ta, Mi}
The finite dual construction $B\mapsto B^\circ$ defines a
contravariant functor on the category of bialgebras that is
adjoint to itself.
\end{Corollary}

\begin{Remark}
If we fix a bialgebra $T$, then the universal bimeasuring
bialgebra construction gives rise to a covariant functor $B(T,\_)$
from the category of commutative algebras to the category of
bialgebras. It is easy to see that the functor preserves
monomorphisms. In particular there is a bialgebra monomorphism
$T^\circ\to B(T,A)$ for any commutative algebra $A$ (arising from
the unit $\eta\colon k\to A$). If the algebra $A$ is augmented,
then the monomorphism is split.
\end{Remark}

\subsection{Bimeasurings over noncommutative algebras}
It makes little sense to discuss bimeasurings when the algebra $A$
is not commutative. A point in case is the following.
\begin{Proposition}\label{comm}
Let $\psi\colon N\ot T\to A$ be a bimeasuring. If either $N$ or
$T$ is a Hopf algebra then $\psi(N\ot T)$ generates a commutative
subalgebra of $A$.
\end{Proposition}
\begin{proof}
Assume $N$ is a Hopf algebra and note that
\begin{eqnarray*}
\psi(n,t)\psi(m,s)&=&\psi(S(m_1),t_1)\psi(m_2n_1,t_2s_1)\psi(S(n_2),s_2)\\
&=& \psi(m,s)\psi(n,t).
\end{eqnarray*}
If $T$ is a Hopf algebra then the argument is symmetric.
\end{proof}

Now suppose that $T$ is a Hopf algebra and that the algebra $A$ is
not commutative. In view of the proposition above, it is clear,
that in case the universal bimeasuring $\theta\colon B(T,A)\ot
T\to A$ can only exist if every bimeasuring from $N\ot T$ to $A$
maps into a fixed commutative subalgebra of $A'$ of $A$. The
proposition below illustrates the fact that the universal
bimeasurings exist in general only if $A$ is abelian.
\begin{Proposition}
The universal bimeasuring bialgebra $\B(k[x],A)$ exists if and only if the
algebra $A$ is commutative.
\end{Proposition}
\begin{proof}
It is sufficient to see that every element of $A$ is in the image
of some bimeasuring $N\ot k[x]\to A$. This is observed by noting
that $\psi_\alpha\colon k[x]\ot k[x]\to A$, given by
$\psi(x^i,x^j)=\delta_{i,j}i!\alpha^i$ is a bimeasuring for all
$\alpha\in A$ ($\delta_{i,j}$ denotes the Kronecker's delta
function).
\end{proof}

\section{Universal cocommutative bimeasuring bialgebras}

\begin{Proposition}
Let $T$ be a bialgebra and $A$ an algebra (not necessarily
commutative). If the universal bimeasuring bialgebra $B(T,A)$
exists, then the universal cocommutative bialgebra $B_c(T,A)$
exists as well and we have the equality $B_c(T,A)=(B(T,A))_c$.
\end{Proposition}
\begin{proof}
Clear.
\end{proof}

Hence if $A$ is a commutative  algebra, then we always have $B_c(T,A)= B(T,A)_c$.
The proposition below sheds some light on the structure of
universal cocommutative bimeasurings.

\begin{Proposition}\label{p1}
Suppose the image of a bimeasuring $\psi\colon N\ot T\to A$ generates a commutative
subalgebra of $A$. If $N$ is cocommutative, then $\psi$ factors through
$T_{ab}$, i.e. there is a unique bimeasuring
$\overline{\psi}\colon N\ot T_{ab}$ such that
$\psi=\overline{\psi}(1\ot\pi)$, where $\pi\colon T\to T_{ab}$ is
the canonical projection.
\end{Proposition}
\begin{proof}
We compute
$$
\psi(n,ts) = \psi(n_1,t)\psi(n_2,s) = \psi(n_2,s)\psi(n_1,t)
=\break \psi(n_1,s)\psi(n_2,t) = \psi(n,st) $$ and conclude the
proof by pointing out that if $\psi(n,t)=0$ for some $t\in T$ and
all $n\in N$, then
$\psi(n,sts')=\psi(n_1,s)\psi(n_2,t)\psi(n_3,s')=0$ for all
$s,s'\in T$.
\end{proof}

\begin{Corollary}
Let $N$ and $T$ be cocommutative bialgebras. If $\psi\colon N\ot
T\to A$ is a bimeasuring with commutative image in $A$, then $\psi$ factors through $N_{ab}\ot
T_{ab}$, i.e. there is a unique bimeasuring $\overline{\psi}\colon
N_{ab}\ot T_{ab}\to A$ such that $\psi=\overline{\psi}(\pi\ot
\pi)$.
\end{Corollary}

\begin{Proposition}
If $T$ is a perfect Hopf algebra (i.e. $T_{ab}=k$) then the universal bimeasuring bialgebra
$B_c(T,A)$ exists for all algebras $A$ and it is equal to the ground field $k$.
\end{Proposition}
\begin{proof}
Apply Lemma \ref{comm} and Proposition \ref{p1}.
\end{proof}

\begin{Proposition}\label{p2}
If $A$ is a commutative algebra and $T$ a cocommutative bialgebra
then the universal bimeasuring bialgebra $B(T,A)$ is commutative.
\end{Proposition}
\begin{proof}
Apply Proposition \ref{p1}.
\end{proof}

It is natural to ask the symmetric question: If $T$ is a
commutative bialgebra, is the universal bimeasuring bialgebra
$B(T,A)$ automatically cocommutative? We conjecture this is not the
case in general, however we can say the following.

\begin{Proposition}
If $A$ is a commutative algebra and $T$ a bialgebra, then
$B_c(T,A) = B_c(T_{ab},A)$
\end{Proposition}
\begin{proof}
Apply Proposition \ref{p1}.
\end{proof}

\begin{Remark}
The proposition above is symmetric to Proposition \ref{p2} in a
sense that the Proposition in question is equivalent saying that
$$B(T_c,A)_{ab}=B(T_c,A).$$
\end{Remark}

\section{Tensor products and universal bimeasurings}
Throughout this section $T$ and $S$ will be bialgebras and $A$ a commutative algebra.
We shall examine how the tensor product $B(T,A)\ot B(S,A)$ of universal bimeasuring bialgebras $B(T,A)$
and $B(S,A)$ is related to the universal bimeasuring bialgebra $B(T\ot S,A)$.
Recall that if the algebra $A$ is the ground field $k$ then we have
$$B(T,k)\ot B(S,k)=T^\circ \ot S^\circ \simeq (T\ot S)^\circ = B(T\ot S,k).$$

We conjecture that, in general, bialgebras $B(T,A)\ot B(S,A)$ and $B(T\ot S,A)$ are not isomorphic.

Since the algebra $A$ is commutative the linear map
$$\psi\colon B(T,A)\ot B(S,A)\ot T\ot S\to A,$$ given by
$$\psi(f\ot g,t\ot s)=\theta_T(f,t)\theta_S(g,s),$$
is a bimeasuring.
Define $$\alpha\colon B(T,A)\ot B(S,A)\to
B(T\ot S,A)$$ to be the unique coalgebra map such that $\psi=\theta_{T\ot S}(\alpha,1)$. Furthermore
let
$\ol{\io_1}\colon B(T\ot S,A)\to B(T,A)$, $\ol{\io_2}\colon B(T\ot S,A)\to B(T,A)$,
$\ol{\pi_1}\colon B(T,A)\to B(T\ot S,A)$ and $\ol{\pi_2}\colon B(S,A)\to B(T\ot S,A)$
be bialgebra maps induced by
$\io_1=1\ot\eta\colon T\to T\ot S$, $\io_2=\eta\ot 1\colon S\to T\ot S$,
$\pi_1=1\ot\ep\colon T\ot S\to T$ and $\pi_2=\ep\ot 1\colon T\ot S\to S$, respectively.
Using Lemma \ref{l1}, it is easy to see that
the following diagram commutes.
$$
\bfig
\putmorphism(0,750)(3,-2)[{B(T,A)}`{B(T\ot S,A)}`{\ol{\pi_1}}]{1000}{1}l
\puthmorphism(0,750)[\phantom{B(T,A)}`{B(T,A)\ot B(S,A)}`{1\ot\eta}]{1000}{1}a
\puthmorphism(1000,750)[\phantom{B(T,A)\ot B(S,A)}`{B(T,A)}`{1\ot \ep}]{1000}{1}a
\putmorphism(2000,750)(-3,-2)[{B(T,A)}`\phantom{B(T\ot S,A)}`{\ol{\io_1}}]{1000}{-1}r
\putvmorphism(1000,750)[\phantom{B(T,A)\ot B(S,A)}`\phantom{B(T\ot S,A)}`{\alpha}]{700}{1}l
\efig
$$
Symmetrically for $B(S,A)$ with $\ol{\io_2}$ and $\ol{\pi_2}$.
Hence the composite map
$$
B(T,A)\ot B(S,A)\stackrel{\alpha}{\to} B(T\ot S,A)\stackrel
{\ol{\io_1}*\ol{\io_2}}{\to} B(T,A)\ot B(S,A)$$
is the identity and therefore $\alpha$ must be an injective mapping. It is easy to see that the restriction
$\alpha\colon B_c(T,A)\ot B_c(S,A)\to B_c(T\ot S,A)$ is an isomorphism.

\section{Cocommutative bimeasurings and Hopf modules}

Throughout this section $T$ and $N$ denote cocommutative Hopf
algebras and $A$ a commutative algebra. Furthermore let $\mu\colon
N\ot T\to N$ and $\nu\colon N\ot T\to T$ be a pair of actions
making $(N,T,\mu,\nu)$ into an abelian matched pair of Hopf
algebras \cite{Mas}. We can then talk about skew bimeasurings
$\psi\colon N\ot T\to A$, that is linear maps satisfying
\begin{eqnarray*}
&&\psi(nm, t)=\psi(n,m_1(t_1))\psi(m_2,t_2),\; \psi(1,t)=\ep(t),\\
&&\psi(n,ts)=\psi(n_1^{t_1},s)\psi(n_2,t),\; \psi(n,1)=\ep(n)
\end{eqnarray*}
(we abbreviate $\mu(n,t)=n(t)$ and $\nu(n,t)=n^t$), or equivalently
\begin{eqnarray*}
&&\psi(ntm)=\psi(nt_1)\psi(t_2m),\; \psi(t)=\ep(t),\\
&&\psi(tns)=\psi(tn_1)\psi(n_2s),\; \psi(n)=\ep(n)
\end{eqnarray*}
(here we identify $nt\in T\bowtie N$ with  $n\ot t\in N\ot T$).
The set $P_{\mu,\nu}(N,T,A)$ of all such maps then becomes an abelian
group under convolution product. It is then easy to observe that
the abelian group of skew bimeasurings is isomorphic to the first
cohomology group of a matched pair (see \cite{GM}, Sections 2.3
and 2.4, for definition and description of cohomology groups
$\mathcal{H}^{*}(N,T,A)$ of the abelian matched pair
$(N,T)=(N,T,\mu,\nu)$ with coefficients in the algebra $A$).

\begin{Proposition}
If $(N,T,\mu,\nu)$ is an abelian matched pair of Hopf algebras, we
have an isomorphism $P_{\mu,\nu}(N,T,A)\simeq
\mathcal{H}^1(N,T,A)$. In particular the abelian group $P(N,T,A)$
of all bimeasurings from $N\ot T$ to $A$ is isomorphic to the
cohomology group $\mathcal{H}^1(N,T,A)$ of the trivial matched
pair $(N,T,1\ot \ep,\ep\ot 1)$.
\end{Proposition}

There is a relation between bimeasurings, Hopf module isomorphisms
and algebras in the category of Hopf modules, which we want to
outline here. A Hopf module $(M,\delta,\mu)$ over a Hopf algebra
$H$ is a $H$-comodule $\delta\colon M\to H\ot M$ together with a
compatible $H$-module structure $\mu\colon H\ot M\to M$, so that
the diagram
$$
\begin{CD}
H\ot M @> \delta_{H\ot M} >> H\ot H\ot M \\
@V \mu VV @V 1\ot\mu VV \\
H @> \delta >> H\ot M
\end{CD}
$$
commutes, i.e. $\delta(hm)=h_1m_{-1}\ot h_2m_0$, where
$\delta_{H\ot M}=(m_H\ot1\ot 1)\tau_{23}(\Delta\ot\delta )$. A morphism of Hopf
modules is just an $H$-linear and $H$-colinear map. The cotensor product
$M\ot^H N$ together with the diagonal action, which restricts from the diagonal
action of $M\ot N$, is a symmetric tensor in the category of Hopf modules
$\Vect^H_H$. The vector space of coinvariants
$$
A=M^{co H}=\mathrm{equ}\left(M\mathop{\rightrightarrows}_{\iota\ot 1}^{\delta} H\ot M\right)
$$
is precisely the image of $\rho=\mu(S\ot 1)\delta\colon M\to M$,
which then has the image factorization $\rho =\kappa\bar\rho\colon
M\to A\to M$, where $\bar\rho\colon M\to A$ is the projection and
$\kappa\colon A\to M$ the inclusion.

\begin{Theorem}\cite{Sw, Mo}\label{t52}
$\theta=(1\ot \bar\rho)\delta\colon M\to H\ot A$ is an isomorphism of Hopf
modules. The functor $(\;)^{co H}\colon \Vect_H^H\to \Vect$ is a tensor
preserving equivalence of categories with inverse
$H\ot\_\colon \Vect\to \Vect_H^H$.
\end{Theorem}

\begin{proof} It is easy to check that $\theta$ is a homomorphism of Hopf
modules and that $\theta\kappa (a)=1\ot a$ for all $a\in A$. It then follows
that $\mu (1\ot\kappa)\theta = id_M$
and $\theta\mu (1\ot\kappa )=id_{H\ot A}$, so that $\theta$ is invertible and
$\theta^{-1}=\mu (1\ot\kappa)$.
\end{proof}

An algebra in $\Vect_H^H$ is a Hopf module $M$ together with Hopf
module maps $\nu\colon H\to M$ and $\nabla\colon  M\ot^H M\to M$
satisfying the usual unitarity and associativity conditions. It
follows that the equivalence described in the preceding theorem
restricts to algebras $(\;)^{co H}\colon \Alg_H^H\to \Alg$.

\begin{Theorem}\label{t53}
If $(M,\delta,\mu)$ is an algebra in $\Alg_H^H$ with algebra of
coinvariants $A$ the the following groups are isomorphic
\begin{enumerate}
\item $\Reg_+(H,A)$, the group of convolution invertible
normalized linear maps $\psi\colon H\to A$, \item $\Aut_A^H(M)$,
the group of $A$-linear $H$-comodule automorphisms $\Phi\colon
M\to M$,
\item the group $\mathcal A$ of $A$-linear action
$\overline{\mu}\colon H\ot M\to M$ such that
$(M,\delta,\overline{\mu})$ is a Hopf module.
\end{enumerate}
\end{Theorem}

\begin{proof} By Theorem \ref{t52} it suffices to consider the $H$-comodule $H\ot A$.
Convolution invertible, normalized linear maps $\psi\colon H\to A$
are in bijective correspondence with $A$-linear $H$-comodule
automorphisms $\phi\colon H\ot A\to H\ot A$, i.e: there is an
isomorphism $\alpha\colon \Reg_+(H,A)\to \Aut_A^H(H\ot A)$ given
by $\alpha(\psi ) = (1\ot m_A)(1\ot\psi\ot 1)(\Delta_H\ot 1)$ and
$\alpha^{-1}(\phi ) = (\ep_H\ot 1)\phi (1\ot\iota_A)$. In
particular, if $\phi =\alpha (\psi )$ then $\phi (h\ot
a)=h_1\ot\psi (h_2)a$.

The $A$-linear $H$-comodule automorphisms $\phi\colon H\ot A\to
H\ot A$ correspond bijectively to $A$-linear actions
$\tilde\mu\colon H\ot H\ot A\to H\ot A$ such that $(H\ot A,
\Delta\ot 1,\tilde\mu )$ is a Hopf module over $H$ with
coinvariants $A$. The bijection is given by the commutative
diagram
$$\begin{CD}
H\ot H\ot A @>\mu >> H\ot A \\
@V1\ot\phi VV @V\phi VV \\
H\ot H\ot A @>\tilde\mu >> H\ot A
\end{CD}$$
i.e: by the isomorphism $\beta\colon \Aut^H_A(H\ot A)\to \mathcal
A$ defined by $\beta (\phi )=\phi\mu (1\ot\phi^{-1})$ and
$\beta^{-1}(\bar\mu )=\bar\mu (1\ot\iota_H\ot 1)$. A tedious, but
straightforward calculation shows that $(H\ot A, \delta ,\bar\mu
)$ is a Hopf module and in fact an algebra in the category of Hopf
modules over $H$. On the other hand, if $\phi =\beta^{-1}(\bar\mu
)=\bar\mu (1\ot\iota_H\ot 1)$ is an $A$-linear $H$-comodule map,
since $\bar\mu $ is an $A$-linear action such that $(h\ot A,\delta, \bar\mu )$
is a Hopf module. By the arguments in the proof of
Theorem 5.2 it follows that $\bar\phi\theta =\id_{H\ot
A}=\theta\bar\phi$. Moreover, $\beta^{-1}\beta (\phi )=\phi\mu
(1\ot\phi^{-1})=\phi$ and $\beta\beta^{-1} (\bar\mu )=\bar\phi\mu
(1\ot\bar\phi^{-1}) =\bar\phi\mu (1\ot\theta
)=\bar\phi\theta\bar\mu =\bar\mu$.
\end{proof}

If $(N,T,\mu ,\nu)$ is a matched pair of cocommutative Hopf algebras with
bismash product $H=T\bowtie N$, then the relation between the the action
$\bar\mu\colon (T\bowtie N)\ot N\ot T\ot A\to N\ot T\ot A$ and the skew bimeasuring
$\psi\colon T\bowtie N\to A$ is given by
\begin{eqnarray*}
\bar\mu (nt\ot m\ot s\ot a)&=&\bar\mu (n\ot\bar\mu (t\ot m\ot s\ot a) \\
&=&\bar\mu (n\ot t_1[m_1]\ot t_2 s\ot\psi (t_3m_2)a) \\
&=&n_1\cdot t_1[m_1]\ot n_2(t_1s_1)\ot\psi (n_3t_2s_2)\psi (t_3m_2)a,
\end{eqnarray*}
where $t[n]=n_2^{S(n_1)(S(t))}=S(S(n)^{S(t)})$.

\begin{Corollary}
If $(N,T,\mu ,\nu )$ is a matched pair of cocommutative Hopf
algebras and $A$ is a commutative algebra then the following
groups are isomorphic:
\begin{enumerate}
\item $\Bimeas (N\ot T,A)$, the group of bimeasurings under
convolution, \item $\Aut^{T\bowtie N}_{N\ot A}(N\ot T\ot
A)\cap\Aut^{T\bowtie N}_{T\ot A}(N\ot T\ot A)$, the group of
$T\bowtie N$-comodule automorphisms which are $N\ot A$-linear as
well as $T\ot A$-linear, \item $\mathcal A$, the group of actions
$\bar\mu\colon (N\bowtie T)\ot (N\ot T\ot A)\to N\ot T\ot A$
diagonal in $N$ as well as in $T$ (i.e: the $N$-action is $N\ot
A$-linear and the $T$-action is $T\ot A$-linear.
\end{enumerate}
\end{Corollary}

\begin{proof}
The result follows directly from Theorem \ref{t53} by a lengthy, routine
computation. We use the identities
\begin{eqnarray*}
t_1[n_1](t_2^{n_2}) &=& \ep(n)t, \\
(t_1[n_1])^{t_2^{n_2}} &=& n\ep(t),
\end{eqnarray*}
connecting distributive law $nt=n_1(t_1)n_2^{t_2}$ and its inverse $tn=t_1[n_1]t_2^{n_2}$,
where $t[n]=S(S(n)^{S(t)})$ and $t^n=S(S(n)(S(t)))$ \cite{GM}.
\end{proof}

\end{document}